%% file: final.tex
\documentclass[12pt]{article}
\usepackage{fullpage,graphicx,psfrag,amsmath,amsfonts,amssymb}
\usepackage[small,bf]{caption}
\usepackage{graphicx}
\usepackage{url,algorithm,amsmath,amsthm}
\usepackage{algorithmicx,algpseudocode}
\usepackage{subcaption}
\input defs.tex

\title{Preconditioning via Diagonal Scaling}
\author{Reza Takapoui \and Hamid Javadi}
\begin{document}
\maketitle

\input intro.tex 
\input equil.tex 

\input example.tex
\input gp.tex
\input param.tex
\input num.tex
\bibliography{final} 
\end{document}

%% file: defs.tex
\newcommand{\BEAS}{\begin{eqnarray*}}
\newcommand{\EEAS}{\end{eqnarray*}}
\newcommand{\BEA}{\begin{eqnarray}}
\newcommand{\EEA}{\end{eqnarray}}
\newcommand{\BEQ}{\begin{equation}}
\newcommand{\EEQ}{\end{equation}}
\newcommand{\BIT}{\begin{itemize}}
\newcommand{\EIT}{\end{itemize}}
\newcommand{\BNUM}{\begin{enumerate}}
\newcommand{\ENUM}{\end{enumerate}}

\newcommand{\BA}{\begin{array}}
\newcommand{\EA}{\end{array}}

\newcommand{\ones}{\mathbf 1}

\newcommand{\reals}{{\mbox{\bf R}}}

\newcommand{\diag}{\mathop{\bf diag}}

{\begin{quote}}{\end{quote}}

\makeatletter
\long\def\@makecaption#1#2{
   \vskip 9pt 
   \begin{small}
   \setbox\@tempboxa\hbox{{\bf #1:} #2}
   \ifdim \wd\@tempboxa > 5.5in
        \begin{center}
        \begin{minipage}[t]{5.5in}
        \addtolength{\baselineskip}{-0.95pt}
        {\bf #1:} #2 \par
        \addtolength{\baselineskip}{0.95pt}
        \end{minipage}
        \end{center}
   \else 
	\hbox to\hsize{\hfil\box\@tempboxa\hfil}  
   \fi
   \end{small}\par
}
\makeatother

\newcounter{oursection}

\newcounter{lecture}

%% file: intro.tex
\section{Introduction}\label{intro}
Interior point methods solve small to medium sized problems to high accuracy in a 
reasonable amount of time. However, for larger problems as well as stochastic 
problems, one needs to use first-order methods such as stochastic gradient descent 
(SGD), the alternating direction method of multipliers (ADMM), and conjugate gradient 
(CG) in order to attain a modest accuracy in a reasonable number of iterations.

The \emph{condition number} of a matrix $A$, denoted by $\kappa(A)$, 
is defined as the ratio of its maximum singular
value to its minimum singular value.
Both theoretical analysis and practical evidence suggest that the precision and 
convergence rate of the first-order methods can depend significantly on the condition 
number of the 
matrices involved in the problems. As an example, the CG
algorithm for solving the linear system $Ax = b$ achieves a faster rate of convergence when
the condition number of $A$ is smaller. Hence, it is desirable
 to decrease the condition number of
matrix $A$ by applying a transformation to it; this process is called 
\emph{preconditioning}.

A special case of preconditioning is called \emph{diagonal scaling}. Here, we are interested 
in finding diagonal matrices $D$ and $E$ to minimize the condition 
number of the matrix $A^{\prime}=DAE$, in order to accelerate first-order methods.
For example, in applying CG to solve $Ax=b$, we can solve $A^\prime \tilde x=Db$ instead,
and recover $x=E^{-1}\tilde x$, while taking advantage of the small condition number of $A^\prime$.

In our numerical experiments and from theoretical analysis,
we have concluded that preconditioning can improve the performance of the first-order methods
in two different ways. First, it can significantly accelerate the linear algebra operations. 
For example in 
\cite{scs}, each step of the algorithm involves running CG which can be done remarkably 
faster if the appropriate preconditioning is applied. The second effect of 
the preconditioning,
which should be distinguished from the first one, is 
decreasing the number of iterations for achieving desired accuracy by following 
different intermediate points in ADMM.

In this report, we first discuss heuristics for diagonal scaling. Next, we motivate 
preconditioning by an example, and then we study preconditioning for 
a specific splitting form in ADMM called \emph{graph projection splitting}. Finally we
examine the performance of our methods by some numerical examples.  

%% file: equil.tex
\section{Matrix equilibration}\label{equil}
Let $A\in\reals^{m\times n}$ be given. \emph{Matrix equilibration} is a heuristic 
method to find diagonal matrices $D\in\reals^{m\times m}$ and 
$E\in\reals^{n\times n}$ to decrease the condition number
of $DAE$. The basic idea in equilibration is to find $D$ and $E$ such 
that all columns of $DAE$ have equal $\ell_p$ norms and all rows 
of $DAE$ have equal $\ell_p$ norms. (The parameter $p\geq1$
can be selected in our algorithm.) We show that this can be formulated as
a convex optimization problem. Consider the problem
\[
\begin{array}{ll}
\mbox{minimize}   & \sum_{i,j}|A_{ij}|^pe^{x_i}e^{y_j}\\
\mbox{subject to}  & \sum_ix_i = 0\\
& \sum_jy_j = 0,
\end{array}
\]
with variables $x\in\reals^m$ and $y\in\reals^n$. The optimality conditions 
for this problem can be shown to be
equivalent to equilibration of $DAE$ for $D = \diag(e^x)$ and $E=\diag(e^y)$.
Although matrix equilibration is a convex problem, it is computationally favorable 
to solve this problem with heuristic iterative methods, rather than using interior
point methods. There are several methods for matrix equilibration,
here we mention a few famous ones. We refer interested readers 
to \cite{brad2011} for a thorough discussion on matrix equilibration.


\BIT
\item \emph{Sinkhorn-Knopp equilibration algorithm}
was originally
designed to convert matrices with nonnegative entries to doubly stochastic
matrices \cite{skequilibration}. However, with a slight modification, it can be used to perform
equilibration in any $\ell_p$ norm.
\item \emph{Ruiz equilibration algorithm} was originally proposed
to equilibrate square matrices \cite{ruizequilibration}. Here, we made a small modification to it so that
$A$ can be rectangular. 
\item \emph{Matrix free algorithms} are the methods
that obtain information about a matrix only through matrix-vector products.
Matrix-free methods are useful when a matrix is represented as an 
operator with no access to the matrix entries \cite{brad2011}.
\EIT
\begin{minipage}[t]{8.2cm}
  \vspace{-15pt}  
  \begin{algorithm}[H]
    \caption{Sinkhorn-Knopp}
\begin{algorithmic}
\State \textbf{Initialize} $d_{1} = \ones_m,  d_{2} =\ones_n$.
\While{$r_{1} > \epsilon_1$ or $r_{2} > \epsilon_2$}
    \State $(d_{1})_i := (Ad_2)^{-1}_i$.
    \State $(d_{2})_j := (A^Td_1)^{-1}_j$.
    \State $B := {\diag(d_{1})}A{\diag(d_{2})}$.
    \State $r_{1} = \frac{\max_i\|B_{i:}\|}
    {\min_i \|B_{i:}\|}, 
    r_{2} = \frac{\max_i\|B_{:j}\|}
    {\min_i \|B_{:j}\|}$.
\EndWhile \\
\Return $D = {\diag(d_1)}, E = {\diag(d_2)}$.
\end{algorithmic}
  \end{algorithm}
\end{minipage}
\begin{minipage}[t]{8.2cm}
  \vspace{-15pt}
  \begin{algorithm}[H]
    \caption{Ruiz}
\begin{algorithmic}
\State \textbf{Initialize} $d_{1} = \ones_m,  d_{2} =\ones_n, B = A$.
\While{$r_{1} > \epsilon_1$ or $r_{2} > \epsilon_2$}
    \State $(d_1)_i := (d_1)_i(\|B_{i:}\|_p)^{-1/2}$.
    \State $(d_2)_j := (d_2)_j (m/n)^{1/2p}(\|B_{:j}\|_p)^{-1/2}$.
    \State $B := {\diag(d_1)}A{\diag(d_2)}$.
      \State $r_{1} = \frac{\max_i\|B_{i:}\|}
    {\min_i \|B_{i:}\|}, 
    r_{2} = \frac{\max_i\|B_{:j}\|}
    {\min_i \|B_{:j}\|}$.
\EndWhile \\
\Return $D = {\diag(d_1)}, E = {\diag(d_2)}$.
\end{algorithmic}
  \end{algorithm}
\end{minipage}
\vspace{5pt}

Although these algorithms work extremely well in practice, there is 
no theoretical guarantee 
that they decrease the condition number. In fact, in our numerical experiments,
we observed that in some cases they might
slightly increase the condition number.


%% file: example.tex
\section{Example}
\label{example}
In this section, we study one simple example to motivate
diagonal scaling for ADMM. Consider the consensus problem 
\begin{equation}
\label{consensus}
\begin{array}{ll}
\mbox{minimize}   & f(x) + g(z)\\
\mbox{subject to}  & x=z,
\end{array}
\end{equation}
with variables $x,z\in\reals^n$, where
$f(x) = (1/2)\|Ax-b\|_2^2$, with $A \in \reals^{m\times n}$, 
$b \in \reals^m$, and $g:\reals^n\rightarrow\reals \cup\{+\infty\}$ is a convex, 
closed, and proper function. This includes several popular problems 
such as lasso, support vector machine, and non-negative least squares.
Instead of solving \eqref{consensus} directly, we can use ADMM to solve the
equivalent problem
\begin{equation}
\begin{array}{ll}
\mbox{minimize}   & f(x) + g(z)\\
\mbox{subject to}  & Fx=Fz,
\label{precconsensus}
\end{array}
\end{equation}
where $F\in\reals^{n\times n}$ is invertible. The augmented Lagrangian for this
problem will be
\begin{equation}
L_F(x,z,y)=f(x)+g(z)+(Fy)^T(x-z)+(1/2)\|x-z\|_{F^2}^2.
\label{augmentedF}
\end{equation}
Defining $Fy$ as the new dual variable, we see that \eqref{augmentedF} is identical
to the augmented Lagrangian for \eqref{consensus} except that a different norm is
used to augment the Lagrangian. 
Notice that taking $F=\sqrt \rho \mathbf I$ is equivalent to scaling the step size in ADMM
by $\rho$. Hence, choosing the
best matrix $F$ (that results in the fastest convergence) 
is at least as difficult as finding the optimal step 
size $\rho$ in \cite{admm:2010}. 
Now, a question of our interest here is the following: by expanding the class 
of matrices $F$ 
from multiples of identity to diagonal matrices, 
how much can we accelerate the ADMM?
To answer this question we note that it can be shown that ADMM achieves 
linear convergence in this case and 
(an upper bound for) the rate of convergence is $\kappa(F(A^TA)^{-1}F^T)$. 
Hence choosing $F$ such that $\kappa(F(A^TA)^{-1}F^T)$ is smaller can
hopefully result in faster convergence for ADMM.

To illustrate this effect, we used ADMM to solve the lasso problem
with and without presence of matrix $F$. We generated a random matrix
$A$ with size $m=7500$ and $n=2500$, and used ADMM to solve
\eqref{precconsensus} in three different cases with the stopping tolerance set to $10^{-4}$.
 First, we performed ADMM directly on \eqref{consensus}
without presence of $F$. It took 7416 iterations 
for the algorithm to converge. For the second case, we took 
$D = \sqrt{\rho^*}\mathbf{I}$, where 
$\rho^* = \sqrt{\lambda_{\min}(A)\lambda_{\max}(A)}$  which maximizes 
the upper bound on the convergence rate of the algorithm ($\rho^*=7.2$ in 
this example). In this case, it took 1019 iterations
until the algorithm converged. 
Finally we used $D$ which equilibrates $(A^TA)^{-1}$. In the last
case it took only 8 iterations for the algorithm to converge. 
This shows the important effect of 
preconditioning on the convergence rate of the ADMM for solving the lasso problem.

%% file: gp.tex
\section{Graph projection splitting}
\label{gp}
Consider the following problem which is in 
canonical \emph{graph form} \cite{graphprojection}
    \begin{equation}
    \begin{array}{ll}
    \mbox{minimize}  & f(y)+g(x) \\
    \mbox{subject to} & Ax=y. \\
    \label{graph}
    \end{array}
    \end{equation}
Graph projection splitting \cite{graphprojection} is a form of ADMM to solve this 
problem serially. The essential idea is to define  $z=(x,y)$ and $\phi(z)=f(y)+g(x)$, and solve
     \begin{equation*}
    \begin{array}{ll}
    \mbox{minimize}  & \phi(z)+ I_{Ax'=b'}(z^\prime) \\
    \mbox{subject to} & z=z^\prime \\
    \label{graph2}
    \end{array}
    \end{equation*}
instead. Here we notice that \eqref{graph} is equivalent to 
\[
    \begin{array}{ll}
    \mbox{minimize}  &  f(D^{-1}\tilde y)+g(E \tilde x)  \\
    \mbox{subject to} & \tilde y = DAE \tilde x, \\
    \end{array} 
\]
where $\tilde x = E^{-1}x$ and $\tilde y = D y$ for diagonal 
$D\in\reals^{m\times m}$ and $E\in\reals^{n\times n}$ with positive diagonal entires.
Using graph projection splitting for this problem is equivalent to running ADMM on
\[
    \begin{array}{ll}
    \mbox{minimize}  & \phi(z)+I_{Ax^\prime=y^\prime}(z^\prime) \\
    \mbox{subject to} & \left[ \begin{array}{cc}
   E^{-1}& 0 \\
   0 & D 
  \end{array} \right]z=\left[ \begin{array}{cc}
   E^{-1}& 0 \\
   0 & D
     \end{array} \right]z^\prime. \\    \end{array}
     \label{graphprjprecond}
\]

The second step of ADMM is to project onto the 
\emph{graph subspace} $G=\left\{(\tilde x,\tilde y)|DAE\tilde x=\tilde y\right\}$ which 
is a linear transformation defined by
\[
\Pi_G(\tilde x,\tilde y)=\left[ \begin{array}{cc}
   I & (DAE)^T \\
   DAE & -I 
 \end{array} \right]^{-1} \left[ \begin{array}{cc}
   I& (DAE)^T \\
   0 & 0 
  \end{array} \right] \left[ \begin{array}{c}
  \tilde x \\ \tilde y \end{array} \right].
\]
In the next section, we will discuss how we can choose $D$ and $E$ to accelerate 
ADMM in graph projection splitting form. Notice that for fixed $DAE$ (even 
when $D$ and $E$ are changing in each step), the projection matrix remains the
same, and hence no additional factorization caching is required in evaluating the projection operator.
We will return to this fact in the next section.

%% file: param.tex
\section{Parameter selection}
\label{param}
In this section we propose methods for choosing matrices $D$, $E$ to 
speed up the graph projection splitting algorithm. We argue that
following the steps below can help us to achieve a desired accuracy in 
smaller number of iterations.

\BIT
\item Using the algorithms discussed in \S\ref{equil}, we choose
$\hat D$, $\hat E$ such that $\hat DA\hat E$ is equilibrated. This will
usually result in reducing the condition number of $A$. As a result, the singular
values of $\hat {D}A\hat E$ will be closer to each other and intuitively
matrix $A$ becomes more isotropic. We notice that defining 
$D=\alpha \hat D$ and $E=\beta \hat E$ 
for any $\alpha, \beta \in \reals_{++}$, the matrix $DAE$ is equilibrated.
We have two degrees of freedom left: $\alpha\beta$ can be chosen to
scale $DAE$ and $\beta/\alpha$ can play the role of step size in ADMM.

\item The product $\alpha \beta$ is chosen to set $\|DAE\|$. 
Let $\gamma = \|\hat DA\hat E\|$, then we will have $\|DAE\|=\alpha
\beta\gamma$.
Since $DAE$ is
likely to have a small condition number, this means that the norms of $y$ and
$x$ are related with a factor close to $\|DAE\|$. This should be chosen
appropriately for every problem.  
In most examples, setting $\|DAE\| = 1$ gives us reasonably good results. Of course,
this depends on the functions $f$ and $g$. We have inspected this effect in the
next section.

\item Note that we can change $\beta/\alpha$ at each step
while $DAE$ is constant. Doing this procedure 
can help us to change $\rho$ adaptively to balance the 
norms of dual and primal residuals. More importantly, this procedure is 
costless and does not require additional factorization caching in direct methods,
also in indirect methods warm start can be used. 
\EIT

%% file: num.tex
\section{Numerical results}
In this section we inspect the effect of parameters explained in \S\ref{param} on the
speed of graph projection splitting algorithm for two different problems.

\subsection{Lasso}
Consider the lasso problem
\[
    \begin{array}{ll}
    \mbox{minimize}  & (1/2)\|Ax-b\|_2^2 + \lambda\|x\|_1, \\
    \label{ll}
    \end{array}
\]
where $A\in\reals^{m\times n}$ and $b\in\reals^m$ and $\lambda\in\reals$ are the 
problem data and $x\in\reals^n$ is the decision variable. 
We can write this problem in form \eqref{graph} with
$f(y)=\|y-b\|_2^2$ and $g(x)=\lambda\|x\|_1$.
We generate instances of the lasso problem with $m=750$ and $n=250$.
We plot the average number of iterations required versus values of $\beta\alpha$ (scaling)
and $\beta/\alpha$ (step size $\rho$) in Figure 1. The relative tolerance was set to $10^{-4}$
as the stopping criterion.

\begin{figure}[ht]
\begin{center}

\psfrag{l1}[B][B]{\raisebox{-1.5ex}{\small$\beta\alpha$}}
\psfrag{l2}[B][B]{\raisebox{1ex}{\small$\beta/\alpha$}}
\psfrag{x1}{\tiny$1/(100\gamma)$}
\psfrag{x2}{\tiny$1/\gamma$}
\psfrag{x3}[L]{\makebox[2ex][r]{\tiny$100/\gamma$}}
\psfrag{y1}[L]{\makebox[1ex][r]{\tiny$\gamma/100$}}
\psfrag{y2}[L]{\makebox[0.25ex][r]{\tiny$\gamma$}}
\psfrag{y3}[T]{\makebox[1ex][r]{\tiny$100\gamma$}}
\includegraphics[width=0.5\textwidth]{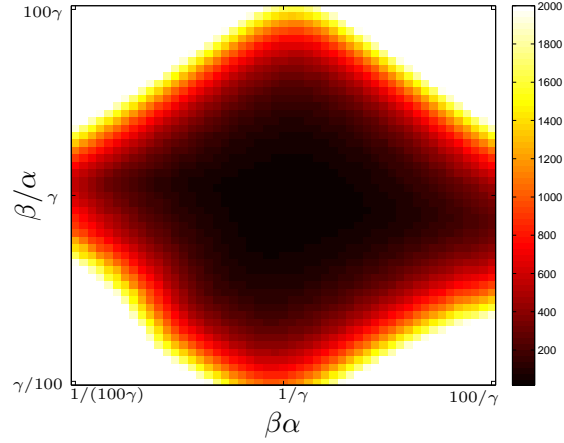}
\end{center}
\caption{\footnotesize Average number of iterations for graph projection splitting on the 
lasso problem}
\end{figure}

\subsection{Linear program (LP)}
Consider the problem
\[
    \begin{array}{ll}
    \mbox{minimize}  & c^Tx \\
    \mbox{subject to} & Ax\preceq b
    \label{lp}
    \end{array}
\]
where $c\in\reals^n$ and $b\in\reals^m$ and $A\in\reals^{m\times n}$
are the problem data and $x\in\reals^n$ is the decision variable. This can be
written in form \eqref{graph} by 
$f(y)=I_{\{y\leq b\}}(y)$ and $g(x)=c^Tx$.
We generate several instances of the problem with $m=750$ and $n=250$. 
We plot the number of iterations required, relative objective suboptimality and
relative constraint error versus values of $\beta\alpha$ (scaling)
and $\beta/\alpha$ (step size $\rho$). The relative tolerance was set to $10^{-4}$
as the stopping criterion.

\begin{figure}[h]
\begin{minipage}{.66\textwidth}
  \begin{subfigure}{\linewidth}
    \centering
    \psfrag{l1}[B][B]{\raisebox{-2ex}{\small$\beta\alpha$}}
	\psfrag{l2}[B][B]{\raisebox{3ex}{\small$\beta/\alpha$}}
\psfrag{x1}[L]{\makebox[6ex][r]{\tiny$50/\gamma$}}
\psfrag{x2}[B][B]{\raisebox{-.25ex}{\tiny$1/\gamma$}}
\psfrag{x3}[L]{\makebox[2ex][r]{\tiny$50/\gamma$}}
\psfrag{y1}[L]{\makebox[1ex][r]{\tiny$1/(50\gamma)$}}
\psfrag{y2}[L]{\makebox[.5ex][r]{\tiny$1/\gamma$}}
\psfrag{y3}[T]{\makebox[1ex][r]{\tiny$50/\gamma$}}
    \includegraphics[width=.6\linewidth]{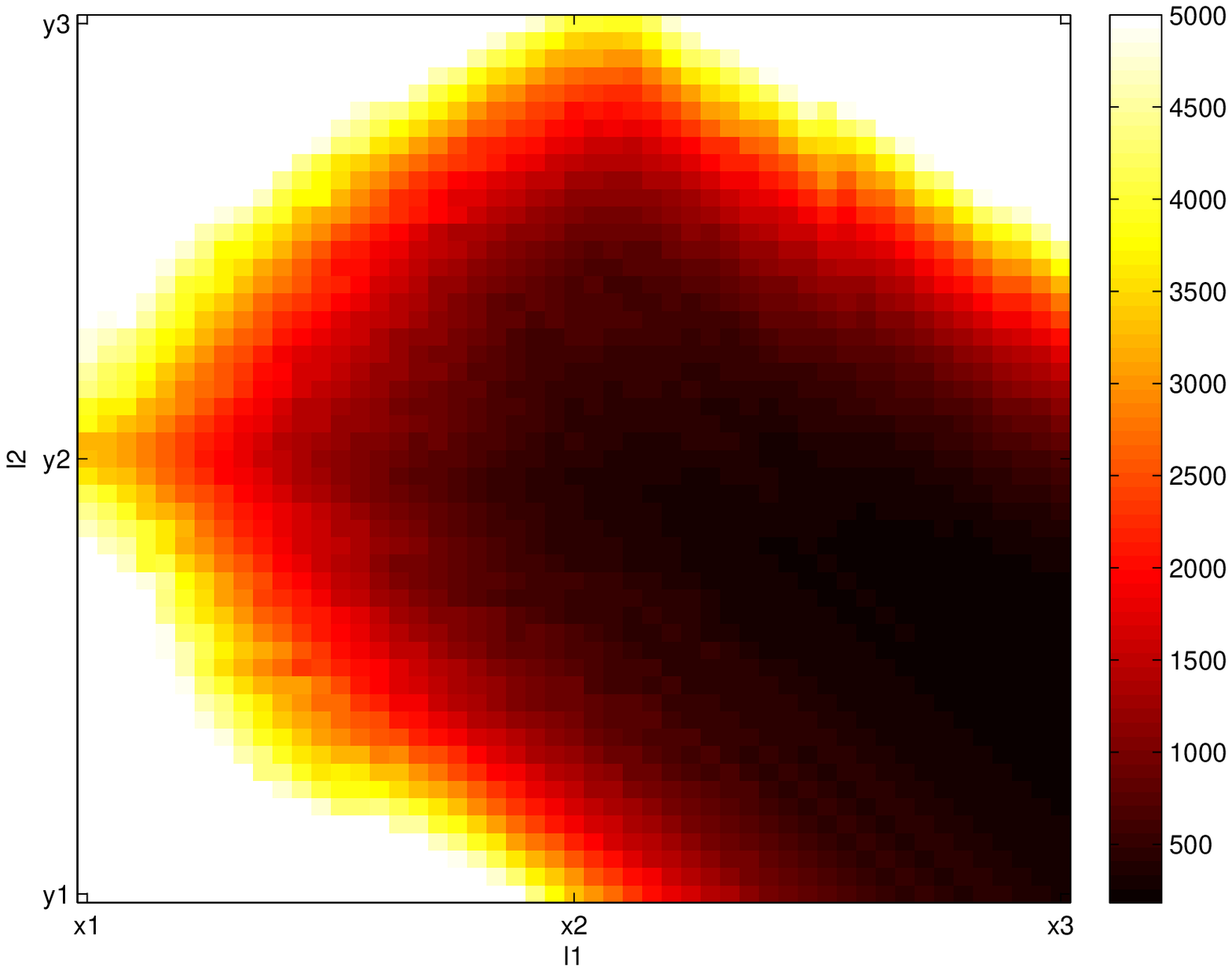}
    \caption{\scriptsize Average number of iterations}
    \label{fig:sub3}
  \end{subfigure}
\end{minipage}%
\begin{minipage}{.33\textwidth}
  \begin{subfigure}{\linewidth}
    \centering
        \psfrag{l1}[B][B]{\raisebox{-1ex}{\tiny$\beta\alpha$}}
	\psfrag{l2}[B][B]{\raisebox{1.5ex}{\tiny$\beta/\alpha$}}
\psfrag{x1}[L]{\makebox[2.8ex][r]{\scalebox{.3}{$50/\gamma$}}}
\psfrag{x2}[B][B]{\raisebox{-.25ex}{\scalebox{.3}{$1/\gamma$}}}
\psfrag{x3}[L]{\makebox[1ex][r]{\scalebox{.3}{$50/\gamma$}}}
\psfrag{y1}[L]{\makebox[0.5ex][r]{\scalebox{.3}{$1/(50\gamma)$}}}
\psfrag{y2}[L]{\makebox[.25ex][r]{\scalebox{.3}{$1/\gamma$}}}
\psfrag{y3}[T]{\makebox[0.3ex][r]{\scalebox{.3}{$50/\gamma$}}}
    \includegraphics[width=.6\linewidth]{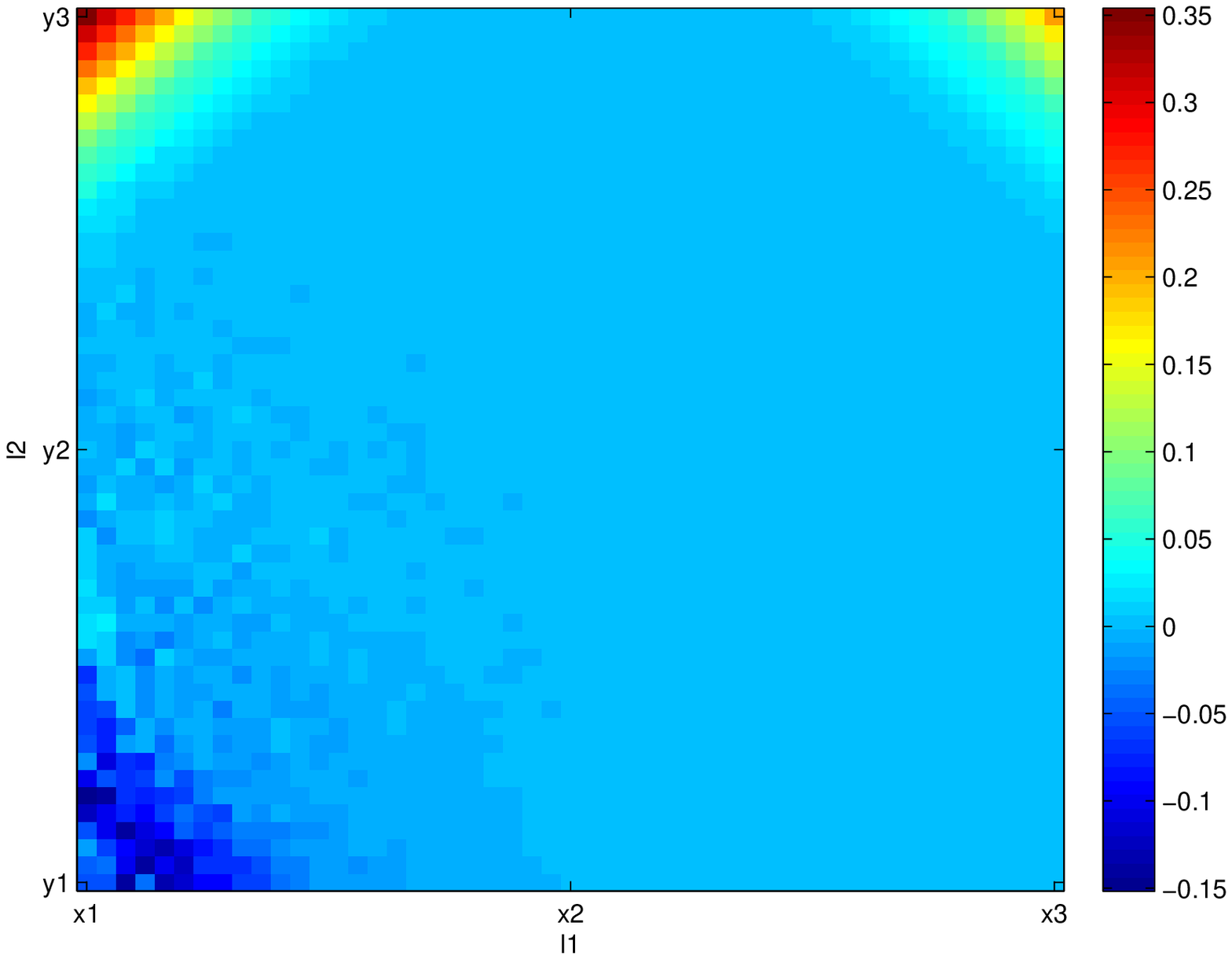}
    \caption{\scriptsize Relative suboptimality}
    \label{fig:sub1}
  \end{subfigure}\\[1ex]
  \begin{subfigure}{\linewidth}
    \centering
        \psfrag{l1}[B][B]{\raisebox{-1ex}{\tiny$\beta\alpha$}}
	\psfrag{l2}[B][B]{\raisebox{1.5ex}{\tiny$\beta/\alpha$}}
\psfrag{x1}[L]{\makebox[2.8ex][r]{\scalebox{.3}{$50/\gamma$}}}
\psfrag{x2}[B][B]{\raisebox{-.25ex}{\scalebox{.3}{$1/\gamma$}}}
\psfrag{x3}[L]{\makebox[1ex][r]{\scalebox{.3}{$50/\gamma$}}}
\psfrag{y1}[L]{\makebox[0.5ex][r]{\scalebox{.3}{$1/(50\gamma)$}}}
\psfrag{y2}[L]{\makebox[.25ex][r]{\scalebox{.3}{$1/\gamma$}}}
\psfrag{y3}[T]{\makebox[0.3ex][r]{\scalebox{.3}{$50/\gamma$}}}
    \includegraphics[width=.6\linewidth]{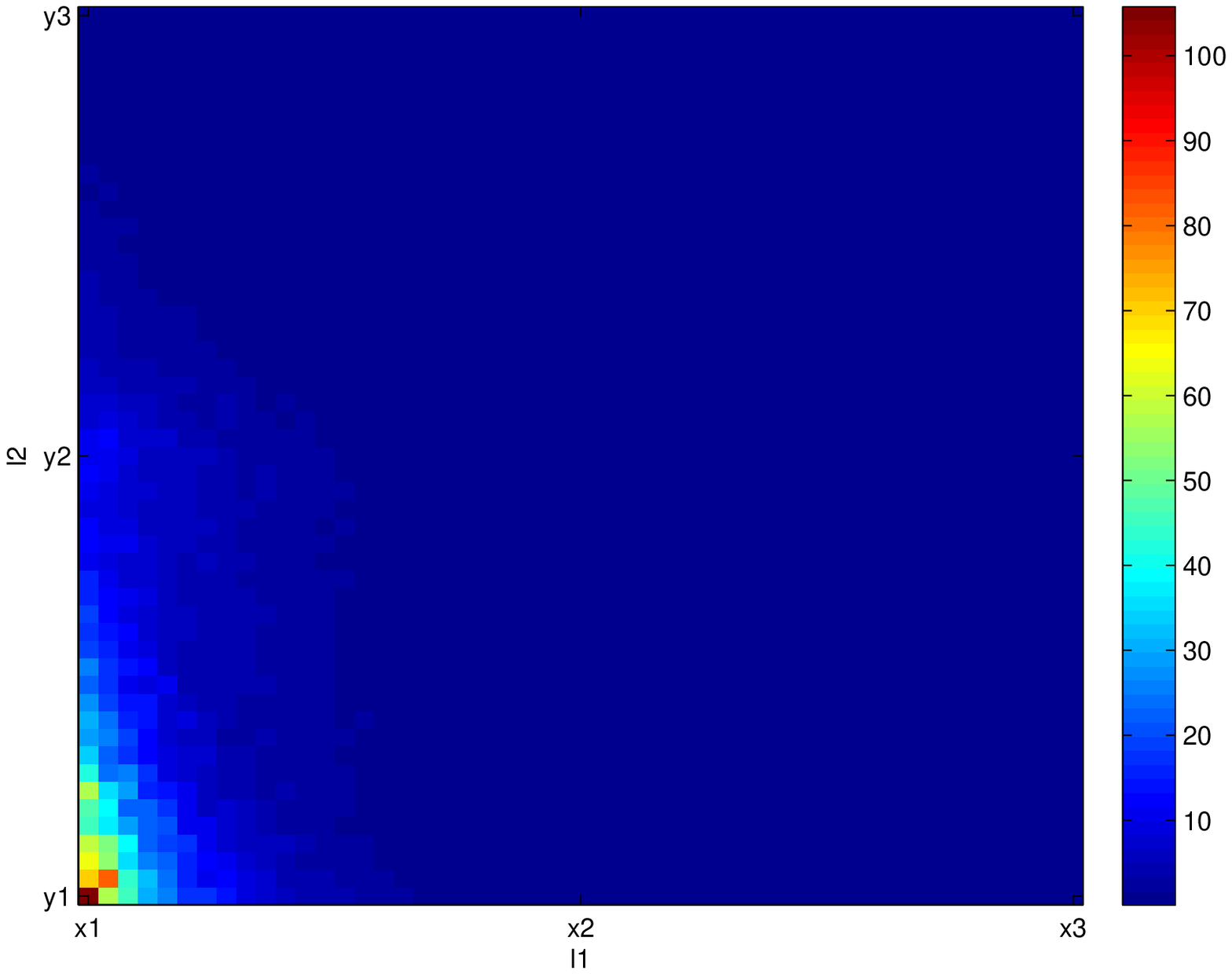}
    \caption{\scriptsize Relative constraint error}
    \label{fig:sub2}
  \end{subfigure}
\end{minipage}%
\caption{\footnotesize Graph projection splitting on LP with inequality constraint}
\end{figure}
We see that after matrix equilibration, the remaining two degrees of freedom
need to be set appropriately to achieve a desired accuracy in smaller number of iterations.
Setting $\alpha$ and $\beta$ such that $\|DAE\|=1$ and choosing appropriate
step size afterwards usually results in a fast convergence.